\documentclass[aos,preprint]{imsart}

\usepackage{xr}
\RequirePackage[OT1]{fontenc}
\usepackage{amsthm,amsmath,cases,bbm,amsfonts,amssymb, dsfont}
\usepackage[authoryear]{natbib}
\RequirePackage[colorlinks,citecolor=blue,urlcolor=blue]{hyperref}
\usepackage{algorithmic,algorithm}

\usepackage[dvipsnames,table,dvipsnames*, svgnames*]{xcolor}

\newcommand{\Ran}{\color{black}}

\usepackage{subcaption}
\usepackage{tocloft}

\usepackage{nccmath, mathtools}

\usepackage{graphicx,comment}
\usepackage{float}
\usepackage{tikz}
\usepackage{mathtools}

\DeclareMathOperator*{\argmin}{\arg\min}

\numberwithin{equation}{section}
\newtheorem{theorem}{Theorem}[section]

\newtheorem{prop}{Proposition}[section]
\theoremstyle{remark}
\newtheorem{example}{Example}[section]
\newtheorem{remark}{Remark}[section]
\newcommand{\beq}{\begin{equation}}
\newcommand{\eeq}{\end{equation}}
\newcommand{\beas}{\begin{eqnarray*}}
\newcommand{\eeas}{\end{eqnarray*}}
\newcommand{\bea}{\begin{eqnarray}}
\newcommand{\eea}{\end{eqnarray}}
\newcommand{\bei}{\begin{itemize}}
\newcommand{\eei}{\end{itemize}}
\newcommand{\ben}{\begin{enumerate}}
\newcommand{\een}{\end{enumerate}}

\newcommand{\FC}{{\mathcal F}}
\newcommand{\F}{{\mathcal F}}
\newcommand{\dd}{\mathrm{d}}
\newcommand{\vv}{m}
\newcommand{\pp}{z}
\newcommand{\CI}{CI}

\newcommand{\E}{\mathbb{E}}

\newcommand{\ih}{ \mathtt{\hat{i}} }
\newcommand{\ti}{\tilde{\mathtt{i}}}
\newcommand{\jh}{\mathtt{\hat{j}}}

\usepackage{etoolbox}
\AfterEndEnvironment{equation}{\restarthintedrel}
\AfterEndEnvironment{align}{\restarthintedrel}
\newcounter{hints}
\renewcommand{\thehints}{\roman{hints}}
\newcommand{\hintedrel}[2][]{%
  \stepcounter{hints}%
  \if\relax\detokenize{#1}\relax\else\csxdef{hint@#1}{\thehints}\fi
  \mathrel{\overset{\textrm{(\thehints)}}{\vphantom{\le}{#2}}}%
}
\newcommand{\restarthintedrel}{\setcounter{hints}{0}}

\begin{document}

\begin{frontmatter}
\title{Estimation and Inference for Minimizer and Minimum of Convex Functions: Optimality, Adaptivity, and Uncertainty Principles}
\runtitle{Minimizer and Minimum of Convex Functions}
 
\begin{aug}
\author{\fnms{T. Tony} \snm{Cai}\ead[label=e1]{tcai@wharton.upenn.edu}},
\author{\fnms{Ran} \snm{Chen}\ead[label=e2]{ran1chen@wharton.upenn.edu}},
\and
\author{\fnms{Yuancheng} \snm{Zhu}\ead[label=e3]{yuancheng.zhu@gmail.com}
\ead[label=u1,url]{URL: http://www-stat.wharton.upenn.edu/$\sim$tcai/}
	\ead[label=u2,url]{URL: http://ran-chen.com/}}
\runauthor{T. T. Cai, R. Chen, and Y. Zhu}
\affiliation{University of Pennsylvania} 
\address{DEPARTMENT OF STATISTICS AND DATA SCIENCE\\
THE WHARTON SCHOOL\\
UNIVERSITY OF PENNSYLVANIA\\
PHILADELPHIA, PENNSYLVANIA 19104\\
\printead{e1}\\
\phantom{E-mail:\ }\printead*{e2}\\
\phantom{E-mail:\ }\printead*{e3}\\
\printead*{u1}\phantom{URL:\ }\\
\printead*{u2}\phantom{URL:\ }
}

\end{aug}

\begin{abstract}
Optimal estimation and inference for both the minimizer and minimum of a convex regression function under the white noise and nonparametric regression models are studied in a non-asymptotic local minimax framework, where the performance of a procedure is evaluated at individual functions.  Fully adaptive and computationally efficient algorithms are proposed and sharp minimax lower bounds are given for both the estimation accuracy and expected length of confidence intervals for the  minimizer and minimum. 

The non-asymptotic local minimax framework brings out new phenomena in simultaneous estimation and inference for the minimizer and minimum.  We establish a novel Uncertainty Principle that provides a fundamental limit on how well the minimizer and minimum can be estimated simultaneously for any convex regression function.  A similar result holds for the expected length of the confidence intervals for the minimizer and minimum.

\end{abstract}



\end{frontmatter}

\section{Introduction}
\label{sec:intro}

Motivated by a range of applications, estimation of and inference for the location and size of the extremum of a nonparametric  regression function has been a longstanding problem in statistics. See, for example,  \cite{kiefer1952stochastic, blum1954multidimensional, chen1988lower}.
The problem has been investigated in different settings. For fixed design, upper bounds for estimating the minimum over various smoothness classes have been obtained  \citep{muller1989adaptive, facer2003nonparametric, shoung2001least}. \citet{belitser2012optimal} establishes the minimax rate of convergence over a given smoothness class for estimating both the minimizer and minimum.
For sequential design, the minimax rate for estimation of the location has been established; see  \cite{chen1996estimation, polyak1990optimal, dippon2003accelerated}. \cite{mokkadem2007companion} introduces a companion for the Kiefer--Wolfowitz--Blum algorithm in sequential design for estimating both the minimizer and minimum.

Another related line of research is the stochastic continuum-armed bandits, which have been used to model online decision problems under uncertainty. Applications include online auctions, web advertising and adaptive routing.  Stochastic continuum-armed bandits can be viewed as aiming to find the maximum of a nonparametric regression function through a sequence of actions.  The objective is to minimize the expected total regret, which requires the trade-off between exploration of new information and exploitation of historical information. See, for example,  \cite{kleinberg2004nearly, Auer2007ImprovedRF, Kleinberg2019}.

In the present paper, we consider optimal estimation and confidence intervals for the minimizer and minimum of convex functions under both the white noise and nonparametric regression models in a non-asymptotic  local minimax framework that evaluates the performance of any procedure at individual functions. This framework provides a much more precise analysis than the conventional minimax theory, which evaluates the performance of the estimators and confidence intervals in the worst case over a large collection of functions. This framework also brings out new phenomena in simultaneous estimation and inference for the minimizer and minimum.

We first focus on the white noise model, which is given by
\beq
\label{white.noise.model}
\dd Y(t) = f(t)\dd t + \varepsilon \dd W(t), \quad 0 \leq t \leq 1,
\eeq
where $W(t)$ is a standard Brownian motion,  and $\varepsilon>0$ is the noise level. The drift function $f$ is assumed to be in $\FC$, the collection of convex functions defined on $[0, 1]$ with a unique minimizer ${\Ran Z(f)} = \argmin_{0\le t\le 1} f(t)$. The minimum value of the function $f$ is denoted by {\Ran $M(f)$}, i.e., ${\Ran M(f)} = \min_{0\le t\le 1} f(t)= f({\Ran Z(f) })$. The goal is to optimally estimate ${\Ran Z(f)}$ and ${\Ran M(f)}$, as well as construct optimal confidence intervals for ${\Ran Z(f)}$ and ${\Ran M(f)}$. Estimation and inference for the minimizer ${\Ran Z(f)}$ and minimum ${\Ran M(f)}$ under the nonparametric regression model will be discussed later in Section \ref{sec:regression}.

\subsection{Function-specific Benchmarks and Uncertainty Principle}
\label{sec:benchmarks}

As the first step toward evaluating the performance of a procedure at individual convex functions in $\FC$,  we define the function-specific benchmarks for estimation of the minimizer and minimum respectively  by
\begin{align}
R_\pp(\varepsilon;f) &= \sup_{g\in\FC} \inf_{\hat Z } \max_{h\in\{f, g\}}\mathbb E_h|\hat Z - Z(h)|,\label{eqn:risk_minimizer} \\
R_\vv(\varepsilon;f) &= \sup_{g\in\FC} \inf_{\hat M} \max_{h\in\{f, g\}}\mathbb E_h|\hat M - M(h)|.\label{eqn:risk_minimimum}
\end{align}

As in \eqref{eqn:risk_minimizer} and \eqref{eqn:risk_minimimum}, we use subscript `$z$' to denote quantities related to the minimizer and `$m$' for the minimum throughout the paper. For any given $f\in \FC$, the benchmarks $R_\pp(\varepsilon;f)$ and $R_\vv(\varepsilon;f)$ quantify the estimation accuracy at $f$ of the minimizer ${\Ran Z(f)}$ and minimum ${\Ran M(f)}$ against the hardest alternative to $f$ within the function class $\FC$.

We show that  $R_\pp(\varepsilon;f)$ and $R_\vv(\varepsilon;f)$ are the right benchmarks for capturing the estimation accuracy at individual functions in $\FC$ and will construct adaptive procedures that simultaneously perform within a constant factor of $R_\pp(\varepsilon;f)$ and $R_\vv(\varepsilon;f)$  for all $f\in \FC$. In addition, it is also shown that any estimator $\hat Z$  for the minimizer that is ``super-efficient" at some $f_0\in \FC$, i.e., it significantly outperforms the benchmark  $R_\pp(\varepsilon;f_0)$,  must pay a penalty at another function $f_1\in \FC$ and thus no procedure can uniformly outperform the benchmark. An analogous result holds for the minimum.

More interestingly, the non-asymptotic local minimax framework enables us to establish a novel Uncertainty Principle for estimating the minimizer and minimum of a convex function. The Uncertainty Principle reveals an intrinsic tension between the task of estimating the minimizer and that of estimating the minimum. That is,  there is a fundamental limit to the estimation accuracy of the minimizer and minimum for all functions in $\FC$ and consequently the minimizer and minimum of a convex function cannot be estimated accurately at the same time.
More specifically, it is shown that
\begin{equation}
R_\pp(\varepsilon;f)\cdot R_\vv(\varepsilon;f)^2  \asymp \varepsilon^2
\end{equation}
for all $f\in \FC$.
Further, on the lower bound side,
\begin{equation}
\label{Estimation-Uncertainty-Principle}
\inf_{f\in \FC} R_\pp(\varepsilon;f)\cdot R_\vv(\varepsilon;f)^2 \geq  {\Phi(-0.5)^3 \over 2} \varepsilon^2,
\end{equation}
where  $\Phi(\cdot)$ is the cumulative distribution function (cdf) of the standard normal distribution.

For confidence intervals with a pre-specified coverage probability, the hardness of the problem is naturally characterized by the expected length.
Let $\mathcal I_{\pp,\alpha}(\F)$  and $\mathcal I_{\vv, \alpha}(\F)$ be, respectively, the collection of confidence intervals for the minimizer ${\Ran Z(f)}$  and the minimum ${\Ran M(f)}$ with guaranteed coverage probability $1 - \alpha$ for all $f\in \F$. Let $L(\CI)$ be the length of a confidence interval $\CI$.
The minimum expected lengths at $f$ of all confidence intervals in $\mathcal I_{\pp,\alpha}(\{f,g\})$ and $\mathcal I_{\vv, \alpha}(\{f,g\})$ with the hardest alternative $g\in \FC$ for $f$ are given by
\begin{align}
L_{\pp,\alpha}(\varepsilon; f) &= \sup_{g\in \FC} \inf_{\CI\in\mathcal I_{\pp,\alpha}(\{f,g\})} \mathbb E_f L(\CI),\label{eqn:ci_minimizer}\\
L_{\vv,\alpha}(\varepsilon; f) &=\sup_{g\in \FC} \inf_{\CI\in\mathcal I_{\vv,\alpha}(\{f,g\})} \mathbb E_f L(\CI).\label{eqn:ci_minimimum}
\end{align}

As in the case of estimation, we will first evaluate these benchmarks for the performance of confidence intervals in terms of the local moduli of continuity and then construct data-driven and computationally efficient confidence interval procedures. Furthermore, we also establish the Uncertainty Principle for the confidence intervals,
\begin{equation}
\label{CI-Uncertainty-Principle}
\inf_{f\in \FC} L_{\pp,\alpha}(\varepsilon; f)\cdot L_{\vv,\alpha}(\varepsilon; f)^2 \geq C_{\alpha}\varepsilon^2.
\end{equation}
where $C_{\alpha}$ is a positive constant depending on $\alpha$ only. The Uncertainty Principle \eqref{CI-Uncertainty-Principle} shows a fundamental limit for the accuracy of simultaneous inference for the minimizer ${\Ran Z(f)}$ and minimum ${\Ran M(f)}$ for any $f\in \FC$.

\subsection{Adaptive Procedures}

Another major step in our analysis is developing data-driven and computationally efficient algorithms for the construction of adaptive estimators and adaptive confidence intervals as well as establishing the optimality of these procedures at each $f\in \FC$.

The key idea behind the construction of the adaptive procedures is to iteratively localize the minimizer by computing the integrals over the relevant subintervals together with a carefully constructed stopping rule.  For estimation of the minimum and minimizer, additional estimation procedures are added after the localization steps. For the construction of the confidence intervals, another important idea is to look back a few steps before the stopping time.

The resulting estimators, ${\Ran \hat{Z}}$ for the minimizer ${\Ran Z(f)} $ and ${\Ran \hat{M}}$ for the minimum ${\Ran M(f)}$, are shown to attain within a constant factor of the benchmarks $R_\pp(\varepsilon;f)$ and $R_\vv(\varepsilon;f)$  simultaneously for all $f\in \FC$,
\[
\mathbb{E}_f |{\Ran \hat{Z}- Z(f)}|  \leq C_z R_\pp(\varepsilon;f)\quad {\rm and} \quad \mathbb{E}_f |{\Ran \hat{M}- M(f)}| \leq C_{m} R_\vv(\varepsilon;f),\\
\]
for some absolute constants $C_z$ and $C_m$ not depending on $f$.
The confidence intervals, $CI_{z,\alpha}$ for the minimizer ${\Ran Z(f)}$ and $CI_{m,\alpha}$ for the minimum ${\Ran M(f)}$, are constructed and shown to be adaptive to individual functions $f\in \FC$, while having guaranteed coverage probability $1-\alpha$. That is, $CI_{z,\alpha}\in \mathcal I_{\pp,\alpha}(\FC)$  and $CI_{ m,\alpha}\in \mathcal I_{\vv, \alpha}(\FC)$ and for all $f\in \FC$,
\beas
\mathbb{E}_f L(CI_{z,\alpha}) &\leq&  C_z(\alpha) L_{\pp,\alpha}(\varepsilon; f) \\
\mathbb{E}_f L(CI_{m,\alpha} ) &\leq& C_m(\alpha) L_{\vv,\alpha}(\varepsilon; f) ,
\eeas
where $C_z(\alpha) $ and $C_m(\alpha)$ are constants depending on $\alpha$ only.

\subsection{Related Literature}

In addition to estimation and inference for the location and size of the extremum of a nonparametric  regression function mentioned at the beginning of this section, the problems considered in the present paper are also connected to nonparametric estimation and inference under shape constraints, which have also been well studied in the literature.

Nonparametric convex regression has been investigated in various settings, ranging from estimation and confidence bands for the whole function \citep{birge1989grenader, guntuboyina2018nonparametric, hengartner1995finite, dumbgen1998new}, to estimation and inference at a fixed point \citep{kiefer1982optimum, cai2013adaptive, cai2015framework, ghosal2017univariate}.  \citet{deng2020inference} established limiting distributions for some local parameters of a convex regression function, including the minimizer based on the convexity-constrained least squares (CLS) estimator and constructed a confidence interval for the minimizer.   As seen in Section \ref{subsec:ci_subopt} of the Supplementary Material \citep{CaiChenZhuSupplement}, this confidence interval is suboptimal in terms of the expected length under the local minimax framework that we introduce later. It is also much more computationally intensive as it requires solving the CLS problem.

In the context of estimating and inferring the value of a convex function at a fixed point, which is a linear functional, the local minimax framework  characterized by the benchmarks \eqref{eqn:risk_minimizer}-\eqref{eqn:risk_minimimum} and \eqref{eqn:ci_minimizer}-\eqref{eqn:ci_minimimum} has been used in  \cite{cai2013adaptive} and  \cite{cai2015framework}.   However, the focus of the present paper is on the minimizer and minimum, which are nonlinear functionals. Due to their nonlinear nature, the analysis is much more challenging than it is for the function value at a fixed point.

Another related line of research is stochastic numerical optimization of convex functions.  \cite{agarwal2011stochastic} studies stochastic convex optimization with bandit feedback and proposes an algorithm that is shown to be nearly minimax optimal.  \cite{zhu2016local} uses the framework introduced in \cite{cai2015framework} to study the local minimax complexity of stochastic convex optimization based on queries to a first-order oracle that produces an unbiased subgradient in a rather restrictive setting.

\subsection{Organization of the Paper}
In Section \ref{sec:lowerbounds}, we analyze individual minimax risks, relating them to appropriate local moduli of continuity and more explicit alternative expressions, and explain the uncertainty principle with a discussion of the connections with the classical minimax framework. Super-efficiency is also considered. In Section \ref{sec:adaptive}, we introduce the adaptive procedures for the white noise model and show that they are optimal. In Section \ref{sec:regression}, we consider the nonparametric regression model. Adaptive procedures are proposed and their optimality is established. Section \ref{sec:discussion} discusses some future directions. For reasons of space, the numerical results and proofs are given in the Supplementary Material \cite{CaiChenZhuSupplement}.

\subsection{Notation}
We finish this section with some notation that will be used in the rest of the paper.
The cdf of the standard normal distribution is denoted by $\Phi$.
For $0<\alpha<1$, $z_{\alpha} = \Phi^{-1}(1-\alpha)$.
For two real numbers $a$ and $b$, $a\wedge b =\min\{a,b\}$, $a\vee b =\max\{a,b\}$.
$\|\cdot\|_2$ denotes the $L_2$ norm (i.e., $\|f\|_2=\sqrt{\int f(x)^2 \dd x}$ ).
For $f\in L_2[0,1]$ and $r>0$, $\mathcal{B}_r(f)=\{g\in L_2[0,1]: \|g-f\|_2\le r\}$ and $\partial{\mathcal{B}_r}(f)=\{g\in L_2[0,1]: \|g-f\|_2= r\}$.

\section{Benchmarks and Uncertainty Principle}
\label{sec:lowerbounds}

In this section, we first introduce the local moduli of continuity and use them to characterize the four benchmarks for estimation and confidence intervals introduced in Section  \ref{sec:benchmarks}, which are summarized in the following table:
\begin{center}
\begin{tabular}{c|cc}
 & Estimation & Inference \\
\hline
Minimizer ${\Ran Z(f)}$ & $R_z(\varepsilon;f)$ & $L_{z,\alpha}(\varepsilon; f)$ \\
Minimum ${\Ran M(f)}$ & $R_m(\varepsilon;f)$ & $L_{m,\alpha}(\varepsilon; f)$. \\
\end{tabular}
\end{center}
We provide alternative expressions for the local moduli of continuity that are easier to evaluate.
The results are used to establish a novel Uncertainty Principle, which shows an intrinsic tension between the estimation/inference accuracy for the minimizer and the minimum for all functions in $\FC$.

\subsection{Local Moduli of Continuity}
\label{sec:local.moduli}

For any given convex function $f\in \FC$, we define the following local moduli of continuity, one for the minimizer, and the other for the minimum,
\bea
\omega_z(\varepsilon; f) &=& \sup\left\{|{\Ran Z(f) - Z(g)}|: \|f - g\|_2 \leq \varepsilon, \, g\in\FC \right\},\label{eqn:modulus_minimizer}\\
\omega_m(\varepsilon; f) &=& \sup\left\{|{\Ran M(f) - M(g)}|: \|f - g\|_2 \leq \varepsilon, \, g\in\FC \right\}, \label{eqn:modulus_minimimum}
\eea
As in the case of  a linear functional,  the local moduli $\omega_z(\varepsilon; f)$ and $\omega_m(\varepsilon; f) $  clearly depend on the function $f$  and can be regarded as an analogue of  the inverse Fisher Information in regular parametric models.

The following theorem characterizes the four benchmarks for estimation and inference in terms of the corresponding local modulus of continuity.

\begin{theorem}
\label{thm:lowerbounds}
Let $0<\alpha < 0.3$. Then
\bea
a_1 \omega_z(\varepsilon; f)  \leq & R_z(\varepsilon; f) & \leq A_1  \omega_z(\varepsilon; f), \\
a_1 \omega_m(\varepsilon; f) \leq &  R_m(\varepsilon; f) &\leq A_1 \omega_m(\varepsilon; f), \\
b_\alpha \omega_z(\varepsilon/3; f)  \leq& L_{z,\alpha}(\varepsilon; f) &\leq B_{\alpha} \omega_z(\varepsilon;f), \\
b_\alpha \omega_m(\varepsilon/3; f) \leq& L_{m,\alpha}(\varepsilon; f) &\leq B_{\alpha} \omega_m(\varepsilon; f),
\eea
where the constants $a_1, A_1, b_\alpha, B_\alpha$ can be taken as $a_1=\Phi(-0.5)\approx 0.309$, $A_1=1.5$,
$b_\alpha=0.6 - 2\alpha$, and $B_{\alpha} = 3(1-2\alpha)z_{\alpha}$.
\end{theorem}

\begin{figure}[htb]
\centerline{\includegraphics[width=4in,height=2in]{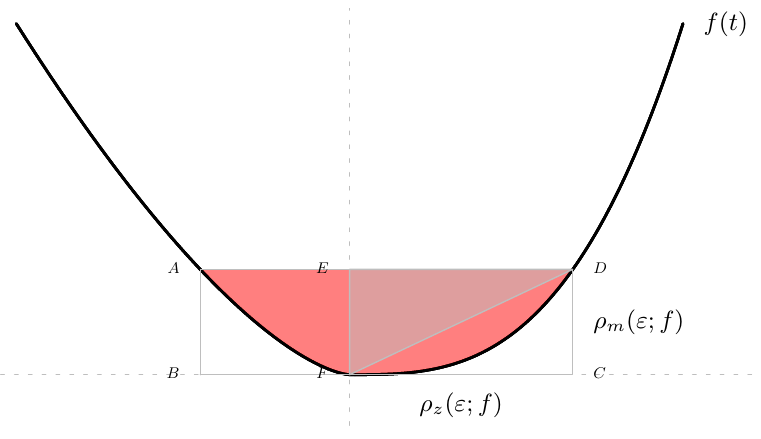}}
\caption{Water filling process.}
\label{fig:water-filling}
\end{figure}

Theorem \ref{thm:lowerbounds} shows that the four benchmarks can be characterized  in terms of the local moduli of continuity. However, these local moduli of continuity are not easy to compute. We now introduce two geometric quantities to facilitate further understanding of these benchmarks. For $f\in \FC$, $u\in \mathbb{R}$ and $\varepsilon > 0$, let $f_u(t) = \max\{f(t), u\}$ and define
\begin{align}
\rho_m(\varepsilon;f) &= \sup\{u - {\Ran M(f)}:\|f-f_u\|_2\leq\varepsilon\},\label{eqn:rho_minimizer}\\
\rho_z(\varepsilon;f) &= \sup\{|t-{\Ran Z(f)}|: f(t)\leq \rho_m(\varepsilon;f) + {\Ran M(f)}, t\in [{\Ran 0,1} ] \}.\label{eqn:rho_minimimum}
\end{align}
Obtaining $\rho_m(\varepsilon;f)$ and $\rho_z(\varepsilon;f)$ can be viewed as a {\it water-filling process.} One adds water into the epigraph defined by the convex function $f$ until the ``volume'' (measured by $\|\cdot\|_2$) is equal to $\varepsilon$.  As illustrated in Figure \ref{fig:water-filling}, $\rho_m(\varepsilon; f)$ measures the depth of the water (CD), and $\rho_z(\varepsilon; f)$ captures the width of the  water surface (FC). $\rho_m(\varepsilon;f)$ and $\rho_z(\varepsilon; f)$ essentially quantify  the flatness of the function $f$ near its minimizer $Z(f)$.

The geometric quantities $\rho_m(\varepsilon;f)$ and $\rho_z(\varepsilon; f)$ defined in \eqref{eqn:rho_minimizer} and \eqref{eqn:rho_minimimum} have the following properties.
\begin{prop}\label{prop:rho}
For $0<c<1$, $f\in \FC$,
\begin{equation}
c \leq { \rho_m(c\varepsilon;f)  \over  \rho_m(\varepsilon;f)} \leq c^{\frac{2}{3}} \quad\text{and}\quad \max\left\{\left({c\over 2}\right)^{\frac{2}{3}}, c \right\} \leq { \rho_z(c\varepsilon;f)  \over  \rho_z(\varepsilon;f)} \leq 1.
\end{equation}
\end{prop}

The following result connects the local moduli of continuity to these two geometric quantities.
\begin{prop}\label{prop:rho_and_omega}
Let $\rho_m(\varepsilon; f)$ and $\rho_z(\varepsilon;f)$ be defined in \eqref{eqn:rho_minimizer} and \eqref{eqn:rho_minimimum}, respectively. Then
\begin{align}
\rho_m(\varepsilon;f) \leq \omega_m(\varepsilon;f) \leq 3\rho_m(\varepsilon;f),\\
\rho_z(\varepsilon;f) \leq \omega_z(\varepsilon;f) \leq 3\rho_z(\varepsilon;f).
\end{align}
\end{prop}

Therefore, through the local moduli of continuity, the hardness of the estimation and inference tasks are tied to the geometry of the convex function near its minimizer.
Note that as the function gets flatter near its minimizer, $\rho_m(\varepsilon; f)$ decreases while $\rho_z(\varepsilon; f)$ increases.
It is useful to calculate $\rho_m(\varepsilon; f)$ and $\rho_z(\varepsilon; f)$ in a concrete example.

\begin{example}
\label{example1}
Consider the function $f(t) = |t-\frac{1}{2}|^k$ where $k\geq 1$ is a constant. We  will calculate $\rho_m(\varepsilon;f)$ and then obtain $\rho_z(\varepsilon;f)$ by first computing $\|f-f_u\|_2^2$ and then setting it to $\varepsilon^2$ to solve for $\rho_m(\varepsilon; f)$.

It is easy to see that in this case
$\|f - f_u\|_2^2 = \frac{4k^2}{(2k+1)(k+1)}\cdot u^{\frac{2k+1}{k}}.$
Setting $\|f - f_u\|_2^2 =\varepsilon^2$ yields
$u=\left(\frac{(2k+1)(k+1)}{4k^2}\right)^{\frac{k}{2k+1}}\varepsilon^{\frac{2k}{2k+1}}.$
Hence,
\[
\rho_m(\varepsilon; f) = \left(\frac{(2k+1)(k+1)}{4k^2}\right)^{\frac{k}{2k+1}}\varepsilon^{\frac{2k}{2k+1}}.
\]
To compute $\rho_z(\varepsilon; f)$, note that
$f^{-1}(u)=\frac{1}{2} \pm u^{1\over k} =\frac{1}{2} \pm  \left(\frac{(2k+1)(k+1)}{4k^2}\right)^{\frac{1}{2k+1}}\varepsilon^{\frac{2}{2k+1}}.$
Hence
\[
\rho_z(\varepsilon; f) = \min\left\{\left(\frac{(2k+1)(k+1)}{4k^2}\right)^{\frac{1}{2k+1}}\varepsilon^{\frac{2}{2k+1}},  \frac{1}{2}\right\}.
\]
Proposition \ref{prop:rho_and_omega} then yields tight bounds for the local moduli of continuity $\omega_m(\varepsilon;f)$ and $\omega_z(\varepsilon;f)$.

\end{example}

\begin{remark}
Note that the results obtained in Example \ref{example1} can be extended to a class of  convex functions. For $f\in \FC$ satisfying $$0<\varliminf\limits_{t \to {\Ran Z(f)}}\frac{f(t)-{\Ran M(f)}}{|t-{\Ran Z(f)}|^k}\le \varlimsup\limits_{t \to {\Ran Z(f)}}\frac{f(t)-{\Ran M(f)}}{|t-{\Ran Z(f)}|^k} < \infty $$ for some $k\ge 1$, it is easy to show that
\[
\omega_m(\varepsilon; f) \sim \varepsilon^{\frac{2k}{2k+1}},\quad \omega_z(\varepsilon; f) \sim \varepsilon^{\frac{2}{2k+1}}, \text{ as } \varepsilon \to 0^+. 
\]
\end{remark}

\subsection{Uncertainty Principle}
\label{sec:uncertainty.principle}

Section \ref{sec:local.moduli} provides a precise characterization of the four benchmarks under the non-asymptotic local minimax framework  in terms of the local moduli of continuity and the geometric quantities $\rho_m(\varepsilon; f)$ and $\rho_z(\varepsilon; f)$. These results yield a novel Uncertainty Principle.

\begin{theorem}[Uncertainty Principle]
\label{thm:uncertainty}
Let $R_z(\varepsilon; f)$, $R_m(\varepsilon; f)$, $L_{z,\alpha}(\varepsilon; f)$, and $L_{m,\alpha}(\varepsilon; f)$ be defined as in \eqref{eqn:risk_minimizer}--\eqref{eqn:ci_minimimum}.  Let $0<\alpha < 0.3$. Then for any $f\in \FC$,
\bea
274\varepsilon^2 >R_z(\varepsilon;f)\cdot R_m(\varepsilon;f)^2 &\geq& \frac{\Phi(-0.5)^3}{2} \varepsilon^2, \label{U.P.z}\\
3^7\cdot(1-2\alpha)^3\varepsilon^2 > L_{z,\alpha}(\varepsilon; f)\cdot L_{m,\alpha}(\varepsilon; f)^2 &\geq& \frac{(0.6-2\alpha)^3}{18} \varepsilon^2. \label{U.P.m}
\eea
\end{theorem}

Note that the bounds in \eqref{U.P.z} and \eqref{U.P.m} are universal for all $f\in \FC$ and show that there is a fundamental limit to the accuracy of estimation and inference for the minimizer and minimum of a convex function.
Our finding here states that the minimizer and the minimum of a convex function cannot be estimated accurately at the same time. This statistical uncertainty principle comes from an intrinsic relationship between the two operators $Z(\cdot)$ and $M(\cdot)$: For any convex function $f\in \FC$ and any $r>0$, there exists $g\in \partial{\mathcal{B}_r}(f)\cap \FC$ such that
\begin{equation}
\label{eq:entangle}
|Z(g)-Z(f)| \cdot |M(g)-M(f)|^2 \ge  \frac{1}{2}\left( {r \over \varepsilon}\right)^2 \cdot \varepsilon^2,
\end{equation}
where $r/\varepsilon=\|(f-g)/\varepsilon\|_2 $ characterizes the probabilistic distance between the two convex functions $f$ and $g$ under the white noise model.

\begin{remark}
\label{rmk:2_generality}
To the best of our knowledge, the uncertainty principles established in this paper are the first of their kind in nonparametric statistics in that they reveal the fundamental tensions between estimation/inference of different quantities. It is shown in the Supplementary Material \citep[Section \ref{subsec:more_on_uncertainty}]{CaiChenZhuSupplement} that similar uncertainty principles also hold for certain subclasses of convex functions.
Note that it is not possible to establish such results using conventional minimax analysis where the performance is measured in a worst-case sense over a large parameter space.
\end{remark}

\subsection{Penalty for Super-efficiency}

We have shown that the estimation benchmarks  $R_z(\varepsilon; f)$ and $R_m(\varepsilon; f)$ defined in \eqref{eqn:risk_minimizer} and \eqref{eqn:risk_minimimum} can be characterized by the local moduli of continuity. Before we show in Section \ref{sec:adaptive} that these benchmarks are indeed achievable by adaptive procedures, we first prove that they cannot be essentially outperformed by any estimator uniformly over $\FC$.
The benchmarks $R_z(\varepsilon;f)$ and $R_m(\varepsilon;f)$ play a role analogous to the information lower bound in classical statistics.

\begin{theorem}[Penalty for super-efficiency]
\label{thm:non_super_efficiency}
For any estimator $\hat{Z}$, if $\mathbb E_{f_0} |\hat{Z}- Z(f_0)|\le \gamma R_z(\varepsilon;f_0)$ for some $f_0\in \FC$ and $\gamma< 0.1$, then there exists $f_1\in \FC$ such that
\beq
\mathbb E_{f_1}( |\hat{Z}- Z(f_1)|)\ge \frac{1}{40} \left( \log{1\over \gamma}\right)^{1/3} R_z(\varepsilon;f_1).
\eeq

Similarly, for any estimator $\hat{M}$, if $\E_{f_0} |\hat{M} - M(f_0)| \le \gamma R_m(\varepsilon;f_0)$  for some $f_0\in \FC$  and $\gamma<0.1$, then there exists $f_1\in \FC$ such that
\beq
\E_{f_1} |\hat{M} - M(f_1)| \ge \frac{1}{8}\left(\log{1\over \gamma}\right)^{1/3} R_m(\varepsilon; f_1).
\eeq
\end{theorem}

\begin{remark}
Theorem \ref{thm:non_super_efficiency} shows that if an estimator of $Z(f)$ or $M(f)$ is super-efficient at some $f_0\in \FC$ in the sense of outperforming the benchmark by a factor of $\gamma$ for some small $\gamma>0$, then it must be sub-efficient at some $f_1\in \FC$ by underperforming the benchmark by at least a factor of $\left( \log{1\over \gamma} \right)^{\frac{1}{3}}$.
\end{remark}

\newcommand{\s}{S_{i_{R}-i_{L},\frac{\alpha}{4}}}
\newcommand{\KT}{\tilde{K}_{\frac{\alpha}{4}}} 
\newcommand{\mjk}{m_{\hat{j}+\KT}}
\newcommand{\sn}{ z_{ \frac{\alpha^{1.84}}{ 468( 2+ z_{\frac{\alpha}{12}}  )^{\frac{2}{3}} }  }}
\newcommand{\tk}{ \max\{4 , 2^{\frac{4}{3}}\times ( 2+z_{\frac{\alpha}{12}})^{\frac{1}{3}} \} }
\newcommand{\scale}{1} 

\section{Adaptive Procedures and Optimality}
\label{sec:adaptive}

We now turn to the construction of data-driven and computationally efficient algorithms for estimation and confidence intervals for the minimizer ${\Ran Z(f)}$ and minimum ${\Ran M(f)}$ under the white noise model. The procedures are shown to be adaptive to each individual function $f\in \FC$ in the sense that they simultaneously achieve, up to a universal constant, the corresponding benchmarks $R_\pp(\varepsilon;f)$, $R_\vv(\varepsilon;f)$, $L_{\pp,\alpha}(\varepsilon; f)$, and $L_{\vv,\alpha}(\varepsilon; f)$ for all $f\in \FC$. These results are much stronger than what can be obtained from a conventional minimax analysis.

\subsection{The Construction}
\label{sec:construction}

There are three main building blocks in the construction of the estimators and confidence intervals: Localization, stopping, and estimation/inference.

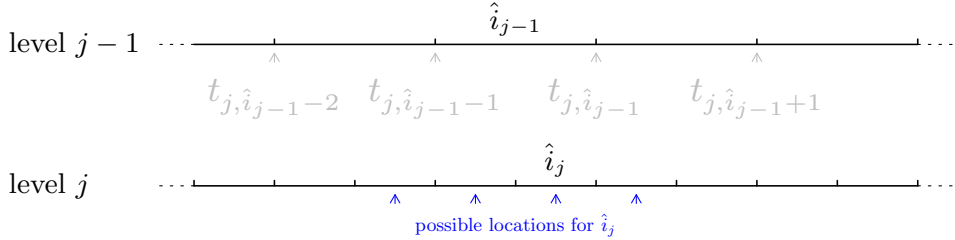
\begin{figure}[H]
\center
\begin{tikzpicture}[x=1pt,y=1pt]
\definecolor{fillColor}{RGB}{255,255,255}
\path[use as bounding box,fill=fillColor,fill opacity=0.00] (0,0) rectangle (361.35,108.41);
\begin{scope}
\path[clip] (  0.00,  0.00) rectangle (361.35,108.41);
\definecolor{drawColor}{RGB}{0,0,0}

\path[draw=drawColor,line width= 0.6pt,line join=round,line cap=round] ( 74.22, 37.47) -- (347.97, 37.47);

\path[draw=drawColor,line width= 0.4pt,dash pattern=on 1pt off 3pt ,line join=round,line cap=round] ( 74.22, 37.47) -- ( 59.01, 37.47);

\path[draw=drawColor,line width= 0.4pt,dash pattern=on 1pt off 3pt ,line join=round,line cap=round] (347.97, 37.47) -- (361.35, 37.47);

\path[draw=drawColor,line width= 0.6pt,line join=round,line cap=round] ( 74.22, 37.47) -- ( 74.22, 39.48);

\path[draw=drawColor,line width= 0.6pt,line join=round,line cap=round] (104.63, 37.47) -- (104.63, 39.48);

\path[draw=drawColor,line width= 0.6pt,line join=round,line cap=round] (135.05, 37.47) -- (135.05, 39.48);

\path[draw=drawColor,line width= 0.6pt,line join=round,line cap=round] (165.47, 37.47) -- (165.47, 39.48);

\path[draw=drawColor,line width= 0.6pt,line join=round,line cap=round] (195.88, 37.47) -- (195.88, 39.48);

\path[draw=drawColor,line width= 0.6pt,line join=round,line cap=round] (226.30, 37.47) -- (226.30, 39.48);

\path[draw=drawColor,line width= 0.6pt,line join=round,line cap=round] (256.72, 37.47) -- (256.72, 39.48);

\path[draw=drawColor,line width= 0.6pt,line join=round,line cap=round] (287.13, 37.47) -- (287.13, 39.48);

\path[draw=drawColor,line width= 0.6pt,line join=round,line cap=round] (317.55, 37.47) -- (317.55, 39.48);

\path[draw=drawColor,line width= 0.6pt,line join=round,line cap=round] (347.97, 37.47) -- (347.97, 39.48);

\node[text=drawColor,anchor=base,inner sep=0pt, outer sep=0pt, scale=  1.20] at (211.09, 43.47) {$\hat i_j$};

\node[text=drawColor,anchor=base west,inner sep=0pt, outer sep=0pt, scale=  1.20] at (  4.18, 34.72) {level $j$};
\definecolor{drawColor}{RGB}{0,0,255}

\node[text=drawColor,anchor=base,inner sep=0pt, outer sep=0pt, scale=  0.80] at (195.88, 20.19) {possible locations for $\hat i_{j}$};

\path[draw=drawColor,line width= 0.4pt,line join=round,line cap=round] (150.26, 30.78) -- (150.26, 34.13);

\path[draw=drawColor,line width= 0.4pt,line join=round,line cap=round] (152.07, 31.00) --
        (150.26, 34.13) --
        (148.45, 31.00);

\path[draw=drawColor,line width= 0.4pt,line join=round,line cap=round] (180.67, 30.78) -- (180.67, 34.13);

\path[draw=drawColor,line width= 0.4pt,line join=round,line cap=round] (182.48, 31.00) --
        (180.67, 34.13) --
        (178.87, 31.00);

\path[draw=drawColor,line width= 0.4pt,line join=round,line cap=round] (211.09, 30.78) -- (211.09, 34.13);

\path[draw=drawColor,line width= 0.4pt,line join=round,line cap=round] (212.90, 31.00) --
        (211.09, 34.13) --
        (209.28, 31.00);

\path[draw=drawColor,line width= 0.4pt,line join=round,line cap=round] (241.51, 30.78) -- (241.51, 34.13);

\path[draw=drawColor,line width= 0.4pt,line join=round,line cap=round] (243.32, 31.00) --
        (241.51, 34.13) --
        (239.70, 31.00);
\definecolor{drawColor}{RGB}{0,0,0}

\path[draw=drawColor,line width= 0.6pt,line join=round,line cap=round] ( 74.22, 91.01) -- (347.97, 91.01);

\path[draw=drawColor,line width= 0.4pt,dash pattern=on 1pt off 3pt ,line join=round,line cap=round] ( 74.22, 91.01) -- ( 59.01, 91.01);

\path[draw=drawColor,line width= 0.4pt,dash pattern=on 1pt off 3pt ,line join=round,line cap=round] (347.97, 91.01) -- (361.35, 91.01);

\path[draw=drawColor,line width= 0.6pt,line join=round,line cap=round] (104.63, 91.01) -- (104.63, 93.01);

\path[draw=drawColor,line width= 0.6pt,line join=round,line cap=round] (165.47, 91.01) -- (165.47, 93.01);

\path[draw=drawColor,line width= 0.6pt,line join=round,line cap=round] (226.30, 91.01) -- (226.30, 93.01);

\path[draw=drawColor,line width= 0.6pt,line join=round,line cap=round] (287.13, 91.01) -- (287.13, 93.01);

\path[draw=drawColor,line width= 0.6pt,line join=round,line cap=round] (347.97, 91.01) -- (347.97, 93.01);

\node[text=drawColor,anchor=base,inner sep=0pt, outer sep=0pt, scale=  1.20] at (195.88, 97.01) {$\hat i_{j-1}$};

\node[text=drawColor,anchor=base west,inner sep=0pt, outer sep=0pt, scale=  1.20] at (  4.18, 88.25) {level $j-1$};
\definecolor{drawColor}{RGB}{190,190,190}

\path[draw=drawColor,line width= 0.4pt,line join=round,line cap=round] (104.63, 84.31) -- (104.63, 87.66);

\path[draw=drawColor,line width= 0.4pt,line join=round,line cap=round] (106.44, 84.53) --
        (104.63, 87.66) --
        (102.83, 84.53);

\path[draw=drawColor,line width= 0.4pt,line join=round,line cap=round] (165.47, 84.31) -- (165.47, 87.66);

\path[draw=drawColor,line width= 0.4pt,line join=round,line cap=round] (167.27, 84.53) --
        (165.47, 87.66) --
        (163.66, 84.53);

\node[text=drawColor,anchor=base,inner sep=0pt, outer sep=0pt, scale=  1.50] at (104.63, 69.71) {$t_{j,\hat i_{j-1}-2}$};

\node[text=drawColor,anchor=base,inner sep=0pt, outer sep=0pt, scale=  1.50] at (165.47, 69.71) {$t_{j,\hat i_{j-1}-1}$};

\path[draw=drawColor,line width= 0.4pt,line join=round,line cap=round] (226.30, 84.31) -- (226.30, 87.66);

\path[draw=drawColor,line width= 0.4pt,line join=round,line cap=round] (228.11, 84.53) --
        (226.30, 87.66) --
        (224.49, 84.53);

\node[text=drawColor,anchor=base,inner sep=0pt, outer sep=0pt, scale=  1.50] at (226.30, 69.71) {$t_{j,\hat i_{j-1}}$};

\path[draw=drawColor,line width= 0.4pt,line join=round,line cap=round] (287.13, 84.31) -- (287.13, 87.66);

\path[draw=drawColor,line width= 0.4pt,line join=round,line cap=round] (288.94, 84.53) --
        (287.13, 87.66) --
        (285.33, 84.53);

\node[text=drawColor,anchor=base,inner sep=0pt, outer sep=0pt, scale=  1.50] at (287.13, 69.71) {$t_{j,\hat i_{j-1}+1}$};
\end{scope}
\end{tikzpicture}
\vspace{-25pt}
\caption{Illustration of the localization step. At level $j$, the middle two intervals are the two subintervals of the selected interval at level $j-1$. One  adjacent interval of the same length on each side is added and the interval at level $j$ is selected among these four intervals.}
\label{fig:localization_strategy}
\end{figure}

In the localization step, we begin with the full interval $[0,1]$. Then, iteratively, we halve the intervals and select one halved interval among a set of halved intervals depending on the interval selected in the previous iteration. This set of halved intervals include the two resulting sub-intervals of the previously selected interval and one neighboring halved interval, when such an interval exists, on both sides. The selection rule is to choose the one with the smallest integral of the white noise process over it. See Figure \ref{fig:localization_strategy} for an illustration of the localization step.

\begin{figure}[H]
\centerline{\includegraphics[scale=1.03]{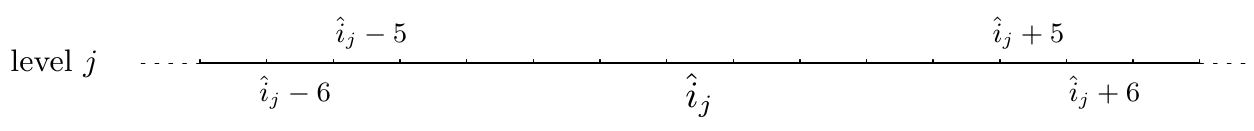}}
\caption{Illustration of the stopping rule.}
\label{fig:stopping_rule}
\end{figure}

The second step of the construction is the stopping rule. The localization step is iterative, so one needs to determine when there is no further gain and stop the iteration. The integral over each selected interval is a random variable and can be viewed as an estimate of the minimum times the length of the interval.
The intuition is that, as the iteration progresses, the bias decreases and the variance increases.  As shown in Figure \ref{fig:stopping_rule}, the basic idea is to use the differences of the integrals over the two neighboring intervals
5 blocks away from the current designated interval, when such intervals exist, on both sides. If either of the differences is smaller than 2 standard deviations, then the iteration stops.

After selecting the final subinterval, the last step in the construction is the estimation/inference for both the minimum and minimizer, which will be described separately later. The detailed construction is given as follows.

\subsubsection{Sample Splitting}

For technical reasons, we split the data into three independent pieces to ensure independence of the data used in the three steps of the construction.
This is done as follows.

Let $B_1(t)$ and $B_2(t)$ be two independent standard Brownian motions, and both be independent of the observed data $Y$. For $t\in [0,1]$, let
\begin{equation}
\begin{split}
Y_l(t)&= Y(t) + \frac{\sqrt{2}}{2} \varepsilon B_1(t) + \frac{\sqrt{6}}{2} \varepsilon B_2(t) , \\
Y_s(t)&= Y(t) + \frac{\sqrt{2}}{2} \varepsilon B_1(t) - \frac{\sqrt{6}}{2} \varepsilon B_2(t) ,\\
Y_e(t)&= Y(t) - \sqrt{2} \varepsilon B_1(t).
\end{split}
\end{equation}
Then $Y_l(\cdot)$, $Y_s(\cdot)$ and $Y_e(\cdot)$ are independent and can be written as
\begin{equation}
\label{eq:sample_splitting_noise}
\begin{split}
dY_l(t) & = f(t)dt + \sqrt{3}\varepsilon dW_1(t) ,\\
dY_s(t) & = f(t)dt + \sqrt{3}\varepsilon dW_2(t) ,\\
dY_e(t) & = f(t)dt + \sqrt{3} \varepsilon dW_3(t) ,
\end{split}
\end{equation}
where $W_1$, $W_2$ and $W_3$ are independent standard Brownian motions.

We now have three independent copies:  $Y_l$ is used for localization, $Y_s$ for stopping, and $Y_e$ for the construction of the final estimator and confidence interval for the minimum.

\begin{remark}
For mere estimation and inference for the minimizer, the copy $Y_e$ is not needed, and it suffices to split into two independent copies with smaller variance. This leads to slightly better performance.
In addition, the noise levels of the three processes $Y_l$, $Y_s$, and $Y_e$ can be different.
For the simplicity and ease of presentation, we split the original sample into three independent and homoskedastic copies for estimation and inference for both the minimizer and minimum.
\end{remark}

\subsubsection{Localization}

For $j = 0, 1, \dots$, and $i = 0, 1, \dots, 2^j$, let
\begin{equation}
\label{local.notation}
m_j = 2^{-j}, \; t_{j,i} =  i \cdot m_j, \;\text{and}\;\, i^*_j=\max\{i: Z(f)\in [t_{j,i-1},t_{j,i}]\}.
\end{equation}
That is, at level $j$ for $j=0, 1, \dots$, the $i^*_j$-th subinterval is the one containing the minimizer $Z(f)$.
For $j=0,1,\dots,$ and $i=1,2,\dots,2^j$, define
\[
X_{j,i} = \int_{t_{j,i-1}}^{t_{j,i}} \dd Y_l(t),
\]
where $Y_l$ is one of the three independent copies constructed above through sample splitting.
For convenience, we define $X_{j,i} =+ \infty$ for $j=0,1,\dots,$ and $i \in \mathbb Z \setminus  \{ 1,2,\dots,2^j\}$.

Let $\hat i_0 = 1$ and for $j = 1, 2, \dots$, let
\begin{align*}
\hat i_j = \argmin_{2\hat i_{j-1}-2 \leq i \leq 2\hat i_{j-1}+1} X_{j,i}.
\end{align*}
Note that given the value of  $\hat i_{j-1}$ at level $j-1$,  in the next iteration the procedure halves the interval $[t_{\hat i_{j-1}-1}, t_{\hat i_{j-1}}]$ into two subintervals and selects  the interval $[t_{\hat i_{j}-1}, t_{\hat i_{j}}]$ at level $j$ from these and their immediate neighboring  subintervals. So $i$ only ranges over 4 possible values at level $j$. See Figure \ref{fig:localization_strategy} for an illustration.

\subsubsection{Stopping Rule}
\label{subsec:stopping_rule}

It is necessary to have a stopping rule to select a final subinterval constructed in the localization iterations. We use another independent copy $Y_s$ constructed in the sample splitting step to devise a stopping rule.
For $j=0,1,\dots,$ and $i=1,2,\dots,2^j$, let
\[
\tilde{X}_{j,i} = \int_{t_{j,i-1}}^{t_{j,i}} \dd Y_s(t).
\]
Again, for convenience, we define $\tilde{X}_{j,i} = +\infty$ for $j=0,1,\dots,$ and $i \in \mathbb Z \setminus  \{ 1,2,\dots,2^j\}$.
Let the statistic $T_j$  be defined as
\[
T_j=\min\{ \tilde{X}_{j,\hat{i}_j + 6} - \tilde{X}_{j, \hat{i}_j+5} ,\;  \tilde{ X}_{j,\hat{i}_j -6} -\tilde{X}_{j,\hat{i}_j-5} \},
\]
where we use the convention $+\infty - x = +\infty $ and $\min\{+\infty, x\}=x$,  for any $-\infty\le x \le \infty$.

The stopping rule is based on the value of $T_j$. It is helpful to provide some intuition before formally defining the stopping rule.
Intuitively, the algorithm should stop at a place where the signal to noise ratio of $T_j$ is small or where the signal is negative.
Let $\sigma_j^2 =6m_j\varepsilon^2$. It is easy to see that, when $\tilde{X}_{j,\hat{i}_j+6} - \tilde{X}_{j,\hat{i}_j+5} < \infty$,
\begin{equation}
 \tilde{X}_{j,\hat{i}_j+6} - \tilde{X}_{j,\hat{i}_j+5} \big| \hat{i}_j \sim N\left(\int_{t_{j,\hat{i}_j+5}}^{ t_{j,\hat{i}_j+6} }\left(f(t+ m_j) -f(t)\right)dt, \; \sigma_j^2 \right).
\end{equation}

Note that the standard deviation $\sigma_j$ decreases at the rate $\frac{1}{\sqrt{2}}$ as $j$ increases. We now turn to the mean of  $\tilde{X}_{j,\hat{i}_j+6} - \tilde{X}_{j,\hat{i}_j+5} | \hat{i}_j$. Recall the notation introduced in \eqref{local.notation}.
It is easy to see that  the algorithm should stop as soon as $\int_{t_{j, \hat{i}_j+5}}^{ t_{j, \hat{i}_j+6} }\left(f(t+ m_j) -f(t)\right)dt$ turns negative, since for any $\hat{i}_j$, if $\int_{t_{j, \hat{i}_j+5}}^{ t_{j, \hat{i}_j+6} }\left(f(t+ m_j) -f(t)\right)dt <0$, then  $|\hat{i}_j - i^*_j|\ge 5$ and consequently $|\hat{i}_{j_1} - i^*_{j_1}|\ge 5$ for any $j_1\ge j$. When $\int_{t_{j, \hat{i}_j+5}}^{ t_{j, \hat{i}_j+6} }\left(f(t+ m_j) -f(t)\right)dt $ is positive, a careful analysis in the proof shows that it shrinks at a rate faster than or equal to ${1\over 4}$ as $j$ increases.
Analogous results hold for  $ \tilde{ X}_{j,\hat{i}_j -6} -\tilde{X}_{j,\hat{i}_j-5}\big| \hat{i}_j $.

Finally, the iterations stop at level $\hat j$ where
\[
\hat{j}= \min\{j: \frac{ T_j}{\sigma_j} \le 2 \}.
\]
The subinterval containing the minimizer $Z(f)$ is localized to be $[t_{\hat{j},\hat{i}_{\hat{j}}-1}, t_{\hat{j},\hat{i}_{\hat{j}}}]$.

\subsubsection{ Estimation and Inference}

After the final subinterval $[t_{\hat{j},\hat{i}_{\hat{j}}-1}, t_{\hat{j},\hat{i}_{\hat{j}}}]$ is obtained, we then use it to construct estimators and confidence intervals for $Z(f)$ and$M(f)$. We begin with the minimizer $Z(f)$.
The estimator of $Z(f)$ is given by the midpoint of the interval $[t_{\hat{j},\hat{i}_{\hat{j}}-1}, t_{\hat{j},\hat{i}_{\hat{j}}}]$, i.e.,
\begin{equation}
\label{hat.Z}
\hat{Z} = \frac{t_{\hat{j},\hat{i}_{\hat{j}}} + t_{\hat{j},\hat{i}_{\hat{j}}-1} }{2}.
\end{equation}

To construct the confidence interval for ${\Ran Z(f)}$, one needs to take a few steps to the left and to the right at level $\hat j$.
Let $K_{\alpha}= \lceil \frac{\log{\alpha}}{\log{\Phi(-2)}} \rceil $ and define
\[
L = \max\{0, \; \hat{i}_{\hat{j}}-12\times 2^{K_{\alpha}}+1\} ,\quad  U = \min\{2^{\hat{j}} , \; \hat{i}_{\hat{j}}+12\times 2^{K_{\alpha}}-2 \} .
\]
The $1-\alpha$ confidence interval for $Z(f)$ is given by
\beq
\label{CI.Z}
CI_{z,\alpha}=[t_{\hat j, L},\; t_{\hat j, U}].
\eeq

For estimation of and confidence interval for the minimum $M(f)$, define
\[
\bar{X}_{j,i} = \int_{t_{j,i-1}}^{t_{j,i}} Y_e(t) dt .
\]
Let
$
{\tilde i }_{\hat{j}}=\hat{i}_{\hat{j}}+2 \left(\mathbbm{1} \{ \tilde{X}_{\hat{j},\hat{i}_{\hat{j}}+6} -\tilde{X}_{j,\hat{i}_{\hat{j}}+5} \le 2\sigma_{\hat j} \} - \mathbbm{1} \{ \tilde{X}_{\hat{j},\hat{i}_{\hat{j}}-6} -\tilde{X}_{j,\hat{i}_{\hat{j}}-5} \le 2\sigma_{\hat j}\}\right)
$
and define the final estimator of the minimum $M(f)$  by
\beq
\label{hat.M}
\hat{M}=\frac{1}{m_{\hat{j}}}\bar{X}_{\hat{j},{\tilde i}_{\hat{j}}}.
\eeq

We now construct confidence interval for ${\Ran M(f)}$.
Recall that $K_{\alpha}= \lceil \frac{\log{\alpha}}{\log{\Phi(-2)}} \rceil$. Compared with the confidence interval for the minimizer, we take four more blocks on each side at the level $(\hat{j}-K_{\frac{\alpha}{4}}-1)_+$.  More specifically, we define
\[
t_{L}=t_{(\hat{j}-K_{\frac{\alpha}{4}}-1)_+ ,\hat{i}_{(\hat{j}-K_{\frac{\alpha}{4}}-1)_+}-5} ,\quad t_{R}=t_{(\hat{j}-K_{\frac{\alpha}{4}}-1)_+ ,\hat{i}_{(\hat{j}-K_{\frac{\alpha}{4}}-1)_+}+4}. 
\]
Set
\beq
\label{tilde.K}
\tilde{K}_{\alpha} = \max\{ 4 , \;  2 + \lceil \log_2{(2+ z_{\alpha/3})}  \rceil \}.
\eeq

Note that at level $\hat{j}+\tilde{K}_{\frac{\alpha}{4}}$, the indices of the intervals with right endpoints $t_L$, $t_R$ respectively are
\[
i_{L} = t_{L} \cdot 2^{ \hat{j}+\tilde{K}_{\frac{\alpha}{4} }}\quad \text{and} \quad  i_{R} = t_{R} \cdot 2^{ \hat{j}+\tilde{K}_{\frac{\alpha}{4} }} .
\]
Note also that $i_{R}-i_{L}=9\times 2^{1+\tilde{K}_{\frac{\alpha}{4}}+K_{\frac{\alpha}{4}}}$, which only depends on $\alpha$.
Define an intermediate estimator of the minimum $M(f)$ by

\[
\hat{f}_1 =\frac{1}{m_{ \hat{j}+\tilde{K}_{\frac{\alpha}{4}} }} \min_{i_{L} < i \le i_{R}} \bar{X}_{ \hat{j}+\tilde{K}_{\frac{\alpha}{4}} , i }.
\]

Let $F_n$ be the cumulative distribution function of $\tilde{v}_n= \max\{v_1,\dots,v_n\}$, where $v_1,\dots,v_n \overset{i.i.d}{\sim} N(0,1)$, and define
\begin{equation}
\label{eq:S_n_beta}
S_{n,\beta}=F_n^{-1}(1-\beta).
\end{equation}
In other words, $S_{n,\beta}$ is the $(1-\beta)$ quantile of the distribution of the maximum of n $i.i.d.$ standard normal variables. Let
\[
f_{lo} = \hat{f}_1 - z_{\alpha/4} \frac{\sqrt{3}\varepsilon}{\sqrt{m_{\hat{j}+\tilde{K}_{\frac{\alpha}{4}}}}} - \frac{\sqrt{3}\varepsilon}{\sqrt{m_{\hat{j}+\tilde{K}_{\frac{\alpha}{4}}}}}, \quad f_{hi} = \hat{f}_1 + S_{i_{R}-i_{L}, \frac{\alpha}{4}} \cdot\frac{\sqrt{3} \varepsilon}{\sqrt{m_{\hat{j}+\tilde{K}_{\frac{\alpha}{4}}}}}.
\]
Then the $(1-\alpha)$ level confidence interval for $M(f)$ is defined as
\beq
\label{CI.M}
CI_{m,\alpha} = [f_{lo},\; f_{hi}].
\eeq

\subsection{Statistical Optimality}

Now we establish the optimality of the adaptive procedures constructed in Section \ref{sec:construction}. The results show that the data-driven estimators and the confidence intervals attain their corresponding local minimax risks/lengths in the sense of staying within a constant multiplier.
We begin with the estimator of the minimizer.
\begin{theorem}[Estimation of Minimizer]
\label{thm:z_estimate}
The estimator $\hat Z$ defined in \eqref{hat.Z} satisfies
\[
\mathbb{E}_f |\hat{Z}-{\Ran Z(f)}| < 53 \rho_z(\varepsilon;f)\le C_z R_z(\varepsilon;f), \quad \text{for all $f\in \FC$},
\]
where  $C_z>0$ is an absolute constant.
\end{theorem}

The following result holds for the confidence interval $CI_{z,\alpha}$.

\begin{theorem}[Confidence Interval for the Minimizer]
\label{thm:z_CI}
Let $0<\alpha<0.3$. The confidence interval $CI_{z,\alpha}$ given in \eqref{CI.Z} is a $(1-\alpha)$ level confidence interval for the minimizer $Z(f)$ and its expected length satisfies
\[
\mathbb{E}_f L(CI_{z,\alpha})\leq (24\times 2^{K_{\alpha}}-3)\times 17.5\times \rho_{z}(\varepsilon;f) \le C_{z,\alpha}  L_{z,\alpha}(\varepsilon;f), \quad \text{for all $f\in \FC$},
\]
where $K_{\alpha}= \lceil \frac{\log{\alpha}}{\log{\Phi(-2)}} \rceil$ and $C_{z,\alpha}$ is a constant depending on $\alpha$ only.
\end{theorem}

Similarly, the estimator and confidence interval for the minimum $M(f)$ are within a constant factor of the benchmarks simultaneously for all $f\in \FC$.
\begin{theorem}[Estimation of Minimum]
\label{thm:m_estimation}
The estimator $\hat M$ defined in \eqref{hat.M} satisfies
\[
\mathbb{E}_f |{\Ran \hat{M}- M(f)}| < 449 \rho_m(\varepsilon;f)\le C_m R_m(\varepsilon;f),  \quad \text{for all $f\in \FC$},
\]
where  $C_m>0$ is an absolute constant.
\end{theorem}

\begin{theorem}[Confidence Interval for the Minimum]
\label{thm:m_CI}
The confidence interval $CI_{m,\alpha}$ given in \eqref{CI.M} is a $(1-\alpha)$ confidence interval for the minimum $M(f)$ and when $0<\alpha<0.3$, its expected length satisfies
\[
\mathbb{E}_f L(CI_{m,\alpha})\leq c_{m,\alpha} \rho_m(\varepsilon;f)\le C_{m,\alpha} L_{m,\alpha}(\varepsilon;f),  \quad \text{for all $f\in \FC$},
\]
where $c_{m,\alpha}$ and $C_{m,\alpha}$ are constants depending on $\alpha$ only.
\end{theorem}

\section{Nonparametric Regression}
\label{sec:regression}

We have so far focused on the white noise model. The procedures and results presented in the previous sections can be extended to nonparametric regression, where we observe
\beq
\label{regression.model}
y_i = f(x_i) + \sigma z_i, \quad i=0,1,2, \ldots,n,
\eeq
with $x_i=\frac{i}{n}$, and $ z_i \overset{i.i.d}{\sim} N(0,1).$ The noise level $\sigma$ is assumed to be known.
The tasks are the same as before: construct optimal estimators and confidence intervals for the minimizer and minimum of $f\in \FC$.

\subsection{Benchmarks and Discretization Errors}
Analogous to the benchmarks for the white noise model defined in Equations (\ref{eqn:risk_minimizer}), (\ref{eqn:risk_minimimum}), (\ref{eqn:ci_minimizer}), (\ref{eqn:ci_minimimum}), we define similar benchmarks for the nonparametric regression model \eqref{regression.model} with $n+1$ equally spaced observations. Denote by $ \mathcal I_{z,\alpha,n}(\mathfrak{F})$ and $ \mathcal I_{m,\alpha,n}(\mathfrak{F})$ respectively the collections of $(1-\alpha)$ level confidence intervals for $Z(f)$ and $M(f)$ on a function class $\mathfrak{F}$ under the regression model \eqref{regression.model} and let
\begin{equation}
\label{eqn:benchmarks}
\begin{split}
\tilde{R}_{z,n}(\sigma;f) =& \sup_{g\in\FC} \inf_{\hat Z }  \max_{h\in\{f, g\}}\mathbb E_h|\hat Z - Z(h)|, \\
\tilde{R}_{m,n}(\sigma;f) =& \sup_{g\in\FC} \inf_{\hat M}  \max_{h\in\{f, g\}}\mathbb E_h|\hat M - M(h)|,\\
\tilde{L}_{z,\alpha,n}(\sigma; f) = & \sup_{g\in \FC} \inf_{\CI\in\mathcal I_{z,\alpha,n}(\{f,g\})} \mathbb E_f L(\CI),\\
\tilde{L}_{m,\alpha,n}(\sigma; f) =& \sup_{g\in \FC} \inf_{\CI\in\mathcal I_{m,\alpha,n}(\{f,g\})} \mathbb E_f L(\CI).
\end{split}
\end{equation}
Compared with the white noise model, estimation and inference for both $Z(f)$ and $M(f)$ incur additional discretization errors, even in the noiseless case. 
See the Supplementary Material \citep[Section \ref{sec:nonparametric_benchmarks_lowerbounds}]{CaiChenZhuSupplement} for further  discussion.

\subsection{Data-driven Procedures}
\label{sec:construction_regression}
Similar to the white noise model, we first split the data into three independent copies and then construct the estimators and confidence intervals for $Z(f)$ and $M(f)$ in three major steps: localization, stopping, and estimation/inference.

\subsubsection{Data Splitting}

Let  $z_{1,0},z_{1,1},\ldots,z_{1,n}, z_{2,0},z_{2,1}, \ldots,z_{2,n}$ be  i.i.d. standard normal random variables, and all be independent of the observed data $\{y_1, ..., y_n\}$. We construct the following three sequences:
\begin{equation}
\begin{split}
&y_{l,i} = y_i + \frac{\sqrt{2}}{2}\sigma z_{1,i} +  \frac{\sqrt{6}}{2}\sigma z_{2,i}, \\
&y_{s,i} = y_i + \frac{\sqrt{2}}{2}\sigma z_{1,i} -  \frac{\sqrt{6}}{2}\sigma z_{2,i} ,\\ 
&y_{e,i} = y_i - \sqrt{2} \sigma z_{1,i} , 
\end{split}
\end{equation}
for $i=0,\ldots,n$. For convenience, let $y_{l,i}=y_{s,i}=y_{e,i}=\infty$ for $i\not\in\{0,1,\ldots,n\}$.
It is easy to see that these random variables are all independent with the same variance $3\sigma^2$ for $i\in\{0,1,\ldots,n\}$. We will use $\{y_{l,i}\}$ for localization, $\{y_{s,i}\}$ for devising the stopping rule, and $\{y_{e,i}\}$ for constructing the final estimation and inference procedures.

Let $J=\lfloor \log_2 (n+1)  \rfloor$.
For $j=0,1,\ldots,J$, $i= 1,2,\ldots, \lfloor \frac{n+1}{2^{J-j-1}} \rfloor $, the $i$-th block at level $j$ consists of $\{x_{(i-1)2^{J-j}},x_{(i-1)2^{J-j}+1},\ldots,x_{i \cdot 2^{J-j}-1}  \}$. Denote the sum of the observations in the $i$-th block at level $j$ for the sequence $u$ ($u=l,s,e$)  as
\[
Y_{j,i,u} = \sum_{k=(i-1)2^{J-j}}^{i\cdot 2^{J-j}-1} y_{u,k}, \text{ for } u=l,s,e.
\]
 Again,  let $ Y_{j,i,u} = +\infty$ when $i \in \mathbbm{Z}\backslash \{1,2,\ldots, \lfloor \frac{n+1}{2^{J-j-1}} \rfloor \}$, for $u=l,s,e$.

\subsubsection{Localization}
We now use $\{y_{l,i},i=0,\ldots,n\}$ to construct a localization procedure.
Let $\ih_0=1$, and for $j=1,2,\ldots, J$, let
\[
\ih_{j} = \argmin_{ \max\{2\ih_{j-1}-2,1\}\le i \le \min\{ 2\ih_{j-1}+1, \lfloor \frac{n+1}{2^{J-j}} \rfloor \}}  Y_{j,i,l}.
\]

This is similar to the localization step in the white noise model. In each iteration, the blocks at the previous level are split into two sub-blocks.  The $i$-th block at level $j-1$ is split into two blocks, the $(2i-1)$-st block and $2i$-th block,  at level $j$. For a given $\ih_{j-1}$, $\ih_{j}$ is the sub-block with the smallest sum among the two sub-blocks of $(j-1)$-level-block $\ih_{j-1}$ and their immediate neighboring sub-blocks.

 \subsubsection{Stopping Rule}
Similar to the  stopping rule for the white noise model, define the statistic $\mathtt{T}_j$ as
\[
\mathtt{T}_j= \min\{Y_{j,\ih_j+6,s}- Y_{j,\ih_j+5,s}  ,  \; Y_{j,\ih_j-6,s}- Y_{j,\ih_j-5,s}  \}.
\]

Let $\tilde{\sigma}_j^2 = 6 \times  2^{J-j} \sigma^2$. It is easy to see that when $Y_{j,\ih_j+6,s}- Y_{j,\ih_j+5,s} <\infty$,
\begin{equation}
Y_{j,\ih_j+6,s}- Y_{j,\ih_j+5,s} | \ih_j \sim N( \sum_{k=(\ih_j+4)2^{J-j}}^{(\ih_j+5)2^{J-j}-1}f(x_{k+2^{J-j}})-f(x_k) , \tilde{\sigma}_j^2).
\end{equation}

Define 
\[
{\check j}=\left\{
\begin{array}{ll}
 \min \{j: \mathtt{T}_j \le 2\tilde{\sigma}_j\}  \quad & \text{if} \;  \{j: \mathtt{T}_j \le 2\tilde{\sigma}_j\}\cap \{ 0,1,2,\ldots,J\} \not=\emptyset \\
 \infty & \text{otherwise}
 \end{array}
 \right. 
\]
and terminate the algorithm at level $\jh=\min\{J,{\check j}\}$. So, either $\mathtt{T}_j$ triggers the stopping rule for some $0\le j \le J$ or  the algorithm reaches the highest possible level $J$.

With the localization strategy and the stopping rule, the final block, the $\ih_{\jh}$-th block at level $\jh$, is given by $\{x_k: (\ih_{\jh}-1) 2^{J-\jh} \le k \le \ih_{\jh} 2^{J-\jh} -1 \}$.

\subsubsection{Estimation and Inference}
\label{sec:4_estimation_inference}

After we have our final block, $\ih_{\jh}$-th block at level $\jh$, we use it to construct estimators and confidence intervals for the minimizer $Z(f)$ and the minimum $M(f)$. We start with the estimation of  $Z(f)$.
The estimator of $Z(f)$ is given as follows:
\begin{equation}
\label{Z.hat.regression}
\hat{Z} =\left\{
\begin{aligned}
-\frac{1}{2n}+\frac{1}{n}( 2^{J-\jh}\ih_{\jh}-2^{J-\jh-1}) ,& \qquad \check{j}<\infty \\
\frac{1}{n}\arg\min_{ \ih_{\jh}-2 \le i \le \ih_{\jh}+2  }y_{e,i-1}   -{1\over n}  ,& \qquad \check{j}=\infty
\end{aligned}
\right.
\end{equation}

To construct the confidence interval for $Z(f)$, we take a few adjacent blocks to the left and right of $\ih_{\jh}$-th block at level $\jh$.
Let
\begin{align*}
\resizebox{\textwidth}{!}{%
$       \mathtt{L} = \max\{0, \ih_{\jh}-12\times 2^{K_{\alpha/2}}+1\}  \quad\text{\rm and}\quad \mathtt{U} = \min\{ \lceil (n+1) 2^{\jh-J} \rceil , \ih_{\jh}+12\times 2^{K_{\alpha/2}}-2 \} $.}
\end{align*}
When $\check{j}<\infty$, let
\begin{align*}
t_{lo} = \frac{2^{J-\jh}}{n} \mathtt{L}-\frac{1}{2n} \quad\text{\rm and}\quad  t_{hi}=\frac{2^{J-\jh}}{n} \mathtt{U}-\frac{1}{2n}.
\end{align*}

When $\check{j}=\infty$, $t_{lo}$ and $t_{hi}$ are calculated by the following Algorithm \ref{alg:t_lo_t_hi}.
Note that $\check{j}=\infty$ means that the procedure is forced to end and  the discretization error can be dominant.

Algorithm \ref{alg:t_lo_t_hi} first iteratively shrinks the original interval $[t_{lo}-{1\over n},t_{hi}+{1\over n}]$ to find the minimizer ${i_m\over n}$ of the function $f$ among the $n+1$ sample points with high probability. In each iteration, the algorithm tests whether the slopes of the segments on both ends are positive or negative. It shrinks the left end with negative slope (on the  left), or shrinks the right end with positive slope (on the right), or stops if no further shrinking is needed on either side.

Note that the minimizer of any convex function with given values at these $n+1$ points is smaller than the intersection of the following two lines:
\begin{equation}
\label{eq:4_lines}
y=f({i_m\over n}) \quad\text{and}\quad y= {f({i_m+2\over n})-f({i_m+1\over n}) \over 1/n }(t-{i_m+1\over n})+ f({i_m+1\over n}).
\end{equation}
Note that these two lines are determined by $f({i_m\over n})$, $f({i_m+1 \over n})$ and $f({i_m+2\over n})$ only. Given the noisy observations at these three points, ${i_m\over n}$, ${i_m+1\over n}$, and ${i_m+2\over n}$, the range of these two lines and the intersection can be inferred, and the right side of the interval can then be shrunk accordingly.

Similar is done for the left side of the confidence interval.
In addition, boundary cases and other complications need to be considered, which are handled in Algorithm \ref{alg:t_lo_t_hi}.

Note that our construction and the theoretical results only rely on convexity. In particular,  the existence of second order derivative is not needed as it is commonly assumed in the literature. This is an important contributing factor to optimality under the non-asymptotic local minimax framework.

\begin{algorithm}[h!]
\caption{Computing $t_{lo}$ and $t_{hi}$ when $\check{j}=\infty$}
\label{alg:t_lo_t_hi}
\begin{algorithmic}
\STATE   $ L \leftarrow \max\{1, \ih_{\jh}-12\times 2^{K_{\alpha/2}} \}-1 $ , $U \leftarrow \min\{n+1, \ih_{\jh}+12\times 2^{K_{\alpha/2}}   \}-1  $, $\alpha_1 \leftarrow \frac{\alpha}{8} $, $\alpha_2=\alpha/24$
\STATE Generate  $z_{3,0},z_{3,1}\ldots,z_{3,n} \overset{i.i.d.}{\sim} N(0,1)$
\STATE $i_l \leftarrow \min\{\{U\}\cup\{i \in [L,U-1]: y_{e,i}+\sqrt{3}\sigma z_{3,i}-(y_{e,i+1}+\sqrt{3}\sigma z_{3,i+1}) \le 2\sqrt{3}\sigma z_{\alpha_1} \}\} $
                 $i_r \leftarrow \max\{\{L-1\} \cup \{i \in [L,U-1]: y_{e,i}+\sqrt{3}\sigma z_{3,i}-(y_{e,i+1}+\sqrt{3}\sigma z_{3,i+1}) \ge -2\sqrt{3}\sigma z_{\alpha_1} \} \}$
\IF{$i_l =U $}
        \IF{$i_l = n$ \text{ and } $ y_{e,n-2}-y_{e,n-1} - \sqrt{3}\sigma (z_{3,n-2} - z_{3,n-1}) + 2\sqrt{6}\sigma z_{\alpha_2} >0 $ } \STATE
        $ t_{hi} \leftarrow 1$\\
        \resizebox{0.92\textwidth}{!}{%
                $t_{lo} \leftarrow \left( (- \frac{ y_{e,n}-y_{e,n-1} - \sqrt{3}\sigma (z_{3,n} - z_{3,n-1})  + 2\sqrt{6}\sigma z_{\alpha_2} }{ n( y_{e,n-2}-y_{e,n-1} - \sqrt{3}\sigma (z_{3,n-2} - z_{3,n-1}) + 2\sqrt{6}\sigma z_{\alpha_2} ) } + \frac{n-1}{n}) \vee \frac{n-1}{n} \right) \wedge \frac{n}{n},$}
        \ELSE \STATE $ t_{lo}=t_{hi}=U/n $
        \ENDIF
\ENDIF
\IF{$i_r =L-1 $}
        \IF{$i_r=-1 $ \text{ and } $ y_{e,2}-y_{e,1} - \sqrt{3}\sigma (z_{3,2} - z_{3,1}) + 2\sqrt{6}\sigma z_{\alpha_2} >0$ }
        \STATE  $t_{hi} \leftarrow \left( ( \frac{ y_{e,0}-y_{e,1} - \sqrt{3}\sigma (z_{3,0} - z_{3,1})  + 2\sqrt{6}\sigma z_{\alpha_2} }{ n( y_{e,2}-y_{e,1} - \sqrt{3}\sigma (z_{3,2} - z_{3,1}) + 2\sqrt{6}\sigma z_{\alpha_2} ) } + \frac{1}{n}) \vee \frac{0}{n} \right) \wedge \frac{1}{n}, t_{lo}=0$
        \ELSE \STATE $ t_{lo}=t_{hi}=0 $
        \ENDIF
\ENDIF
\IF{$(i_l-U)(i_r-L+1)\not=0$}
   \STATE        $i_{lo} \leftarrow (i_l-1) \vee L  $, $i_{hi} \leftarrow (i_r+2) \wedge U  $
        \IF{ $ i_{hi} - i_{lo}\ge 3$ or $(i_{hi}-n)i_{lo}=0 $  }
        \STATE $t_{lo}=i_{lo}/n $, $t_{hi}=i_{hi}/n$
                        \ELSIF{ $y_{e,i_{hi}+1}-y_{e,i_{hi}}- \sqrt{3}\sigma (z_{3,i_{hi}+1} - z_{3,i_{hi}}) \le -2\sqrt{6}\sigma z_{\alpha_2} $ \text{ or } $y_{e,i_{lo}-1}-y_{e,i_{lo}}- \sqrt{3}\sigma (z_{3,i_{lo}-1} - z_{3,i_{lo}}) \le -2\sqrt{6}\sigma z_{\alpha_2} $ }
                                \STATE $t_{lo}=t_{hi} =(i_{hi}+i_{lo})/2n$
                                \ELSE
                                        \STATE
                                        \resizebox{0.92\textwidth}{!}{%
                                                $t_{hi} \leftarrow \left( ( \frac{ y_{e,i_{hi}-1}-y_{e,i_{hi}} - \sqrt{3}\sigma (z_{3,i_{hi}-1} - z_{3,i_{hi}})  + 2\sqrt{6}\sigma z_{\alpha_2} }{ n( y_{e,i_{hi}+1}-y_{e,i_{hi}} - \sqrt{3}\sigma (z_{3,i_{hi}+1} - z_{3,i_{hi}}) + 2\sqrt{6}\sigma z_{\alpha_2} ) } + \frac{i_{hi}}{n}) \vee \frac{i_{hi}-1}{n} \right) \wedge \frac{i_{hi}}{n}$
                                                }\\
                                         \resizebox{0.92\textwidth}{!}{%
                                         $t_{lo} \leftarrow \left( (- \frac{ y_{e,i_{lo}+1}-y_{e,i_{lo}} - \sqrt{3}\sigma (z_{3,i_{lo}+1} - z_{3,i_{lo}})  + 2\sqrt{6}\sigma z_{\alpha_2} }{ n( y_{e,i_{lo}-1}-y_{e,i_{lo}} - \sqrt{3}\sigma (z_{3,i_{lo}-1} - z_{3,i_{lo}}) + 2\sqrt{6}\sigma z_{\alpha_2} ) } + \frac{i_{lo}}{n}) \vee \frac{i_{lo}}{n} \right) \wedge \frac{i_{lo}+1}{n}$
                                         }
        \ENDIF
\ENDIF
\end{algorithmic}
\end{algorithm}

\smallskip
The $(1-\alpha)$-level confidence interval for the minimizer $Z(f)$ is given by
\begin{equation}
\label{Z.CI.regression}
\mathtt{CI}_{z,\alpha}=[t_{lo}\wedge t_{hi} ,t_{hi}] 
\end{equation}

We now construct the estimator and confidence interval for the minimum $M(f)$. Let $\Delta=\mathbbm{1}\{ Y_{\jh,\ih_{\jh}+6,s} -Y_{\jh,\ih_{\jh}+5,s} \le 2 \sqrt{6}\sigma\sqrt{2^{J-\jh}} \} -   \mathbbm{1}\{ Y_{\jh,\ih_{\jh}-6,s} -Y_{\jh,\ih_{\jh}-5,s} \le 2 \sqrt{6}\sigma\sqrt{2^{J-\jh}} \} $ and define 
\begin{equation}
\ti_{\jh}= \left\{
\begin{array}{ll}
 \ih_{\jh}+ 2 \Delta &\quad  \text{if $ {\check{j}} < \infty$}\\ 
\argmin_{\ih_{\jh}-2 \le i\le \ih_{\jh}+2 }  y_{e,i-1} &\quad  \text{if $ {\check{j}} = \infty$}
\end{array}
\right. .
\end{equation}

The estimator of $M(f)$ is then given by the average of the observations of the copy for estimation and inference in the $\ti_{\jh}$-th block at level $\jh$, 
\begin{equation}
\label{M.hat.regression}
\hat{M} = \frac{1}{2^{J-\jh}}  Y_{\jh,\ti_{\jh},e}.
\end{equation}

To construct the confidence interval for $M(f)$, we specify two levels $j_s$ and $j_l$, with
\[
j_s = \max\{0, \jh-K_{\frac{\alpha}{4}}-1 \} \quad \text{and}\quad  j_l = \min\{ J, \jh+ \tilde{K}_{\frac{\alpha}{4}} \},
\]
where $\tilde{K}_{\frac{\alpha}{4}} $ is defined as in Equation \eqref{tilde.K}. It will be shown that at level $j_s$,  $Z(f)$ is within four blocks of the chosen block with probability at least $1-\frac{\alpha}{4}$, and at level $j_l$,  with probability at least $1-\frac{\alpha}{4}$,  the length of the block is no larger than $\rho_z(\frac{\sigma}{\sqrt{n}};f)$.
Define 
\[
I_{lo} =\max\{1 , 2^{j_l-j_s}(\ih_{j_s}-5) \},\quad  I_{hi}=\min\{ 2^{j_l-j_s}(\ih_{j_s}+4)+1, \lceil \frac{n+1}{2^{J-j_l}} \rceil \}.
\]
It can be shown that $Z(f)$ lies with high probability in the interval \\ \mbox{ $[ \frac{2^{J-j_l}(I_{lo}-1)}{n} ,  \frac{2^{J-j_l}I_{hi}-1}{n}]\cap[0,1]$}.
Define an intermediate estimator for $M(f)$ by 
\[
\mathtt{\hat{f}}_1 = \min_{ I_{lo} \le i \le I_{hi}}  \frac{1}{2^{J-j_l}} Y_{j_l,i,e}.
\]
Let 
\[
\mathtt{f}_{hi}= \mathtt{\hat{f}}_1 + S_{ I_{hi}-I_{lo}+1 ,\frac{\alpha}{4}}\frac{\sqrt{3}\sigma}{\sqrt{ 2^{ J-j_l } }}
\]
 where $S_{n,\beta}$ is defined in Equation \eqref{eq:S_n_beta} in Section \ref{sec:adaptive}.  This is the upper limit of the confidence interval, now we define the lower limit $\mathtt{f}_{lo}$.

When $\jh+\tilde{K}_{\frac{\alpha}{4}}\le J$, let 
\begin{align*}
\mathtt{f}_{lo}= \mathtt{\hat{f}}_1 -(z_{\alpha/4}+1)\frac{\sqrt{3}\sigma}{\sqrt{2^{J-j_l }}}.
\end{align*}

When $\jh+\tilde{K}_{\frac{\alpha}{4}}> J$, we compute $\mathtt{f}_{lo}$ by Algorithm \ref{alg:f_lo}, which is based on the geometric property of the convex function $f$ that for any $1 \le k\le n-2$,
\begin{equation}
        \medmath{
\begin{multlined}
\label{eq:4_lines2}
\inf_{t\in [x_k, x_{k+1}]} f(t) 
\ge  \inf_{t\in[{k\over n},{k+1\over n}]}\max  \Big\{ {f(x_{k})-f(x_{k-1}) \over 1/n}(t-x_{k})+f(x_{k}),\\
{f(x_{k+2})-f(x_{k+1}) \over 1/n}(t-x_{k+1})+f(x_{k+1}) \Big\} .  
\end{multlined}
        }
\end{equation}

\begin{algorithm}[htb]
\caption{ Computing $\mathtt{f}_{lo}$ when $\jh+\tilde{K}_{\frac{\alpha}{4}}> J$}
\label{alg:f_lo}
\begin{algorithmic}
\STATE  $H \leftarrow S_{I_{hi}-I_{lo}+3, \frac{\alpha}{8}}\sqrt{3}\sigma $, $k_l \leftarrow I_{lo}-1$, $k_r \leftarrow I_{hi}-2$
\IF{$I_{lo}=1$}
        \STATE $v_{r,0}(t) \leftarrow { y_{e,2}-y_{e,1}+2H \over 1/n } (t-1/n)+y_{e,1}-H$,$h(0)\leftarrow \min_{t\in[0,1/n]} v_{r,0}(t)$,  $k_l \leftarrow I_{lo}$
        \ENDIF
\IF{$I_{hi} -1 = n$}
        \STATE  $v_{l,n-1}(t) \leftarrow { y_{e,n-1}-y_{e,n-2} -2H \over 1/n }(t-{n-1 \over n})+y_{e,n-1}-H   $, $h(n-1)= \min_{t\in[{n-1\over n},1]}v_{l,n-1} (t)$, $k_r \leftarrow I_{hi}-3$
        \ENDIF
\FOR {$i=k_l,\ldots,k_r$}
        \STATE  Define two linear functions:
        \STATE \resizebox{0.95\columnwidth}{!}{%
                $v_{l,i}(t)= \frac{y_{e,i} -y_{e,i-1}-2H}{1/n} (t-x_i)+ y_{e,i}-H  \text{, } v_{r,i}=  \frac{y_{e,i+2} -y_{e,i+1}+2H}{1/n}(t-x_{i+1})+y_{e,i+1}-H $}
                         \STATE $h(i)=\min_{t\in[x_i,x_{i+1}]}\max\{v_{l,i}(t),v_{r,i}(t)\} $
        \ENDFOR
\STATE $\mathtt{f}_{lo}  \leftarrow \min\{h(i): I_{lo}-1 \le i \le I_{hi}-2\} \wedge \mathtt{ f}_{hi}$
\end{algorithmic}
\end{algorithm}

Note that $h(i)$ in Algorithm \ref{alg:f_lo} is derived from one or two linear functions, so given the relationship of the function values at two end points of the corresponding interval, it has an explicit form. Hence the procedure is still computationally efficient.

The $(1-\alpha)$-level confidence interval for the minimum $M(f)$ is given by 
\begin{equation}
\label{M.CI.regression}
\mathtt{CI}_{m,\alpha} = [\mathtt{f}_{lo}, \; \mathtt{f}_{hi}].
\end{equation}

\begin{remark}
\label{rmk:4_comparison}
As mentioned in the introduction, \cite{agarwal2011stochastic} proposes an algorithm for stochastic convex optimization with bandit feedback. While both our procedures and the method in \citet{agarwal2011stochastic} include an ingredient trying to localize the minimizer through shrinking intervals by exploiting the convexity of the underlying function, the two methods are essentially different due to the significant differences in both the designs and loss functions. The goal of exploiting convexity in \cite{agarwal2011stochastic} is mainly to determine the direction of shrinking their intervals, while ours is mainly for deciding when to stop and what to do after stopping.
\end{remark}

\subsection{Statistical Optimality}
\label{sec:regression_optimality}

Now we establish the optimality of the adaptive procedures constructed in Section \ref{sec:construction_regression}. The regression model is similar to the white noise model, but with additional discretization errors. The results show that our data-driven procedures are simultaneously optimal (up to a constant factor) for all $f\in \FC$. We begin with the estimator of the minimizer.

\begin{theorem}[Estimation of the Minimizer]
\label{thm:r_z_estimation}
The estimator $\hat{Z}$ of the minimizer $Z(f)$ defined in (\ref{Z.hat.regression}) satisfies
\beq
\label{eq:reg-minimizer}
\mathbb{E}_f |\hat{Z}-Z(f)|   \le C_1\tilde{R}_{z,n}(\sigma;f),  \quad \text{for all $f\in \FC$},
\eeq
where $C_1>0$ is an absolute constant.
\end{theorem}

The following result holds for the confidence interval $\mathtt{CI}_{z,\alpha}$ of $Z(f)$.
\begin{theorem}
\label{thm:r_z_CI}
Let $0<\alpha<0.3$. The confidence interval $\mathtt{CI}_{z,\alpha}$ given in (\ref{Z.CI.regression}) is a $(1-\alpha)$-level confidence interval for the minimizer $Z(f)$ and its expected length satisfies
$$\mathbb E_f  L(\mathtt{CI}_{z,\alpha})  \le C_{2,\alpha}\tilde{L}_{z,\alpha,n}(\sigma; f), \quad \text{for all $f\in \FC$},$$
where $C_{2,\alpha}$ is a constant depending on $\alpha$ only.
\end{theorem}

Similarly, the estimator and confidence interval for the minimum $M(f)$ are within a constant factor of the benchmarks simultaneously for all $f\in \FC$.

\begin{theorem}[estimation for the minimum]
\label{thm:r_m_estimation}
The estimator $\hat{M}$ defined in (\ref{M.hat.regression}) satisfies
$$ \mathbb E_f | \hat{M} - M(f) |  \le C_3 \tilde{R}_{m,n}(\sigma;f), \quad \text{for all $f\in \FC$},$$
where $C_3$ is an absolute constant.
\end{theorem}

\begin{theorem}
\label{thm:r_m_CI}
Let $0<\alpha<0.3$. The confidence interval $\mathtt{CI}_{m,\alpha}$ given in (\ref{M.CI.regression}) is a $(1-\alpha)$-level confidence interval and its expected length satisfies
$$ \mathbb E_f L(\mathtt{CI}_{m,\alpha}) \le C_{4,\alpha} \tilde{L}_{m,\alpha,n}(\sigma; f), \quad \text{for all $f\in \FC$},$$
where $C_{4,\alpha}$ is a constant depending only on $\alpha$.
\end{theorem}

In addition to statistical optimality, the proposed algorithms are computationally fast. We conducted numerical experiments and the results are in the Supplementary Material \citep[Section  \ref{sec:simulation}]{CaiChenZhuSupplement}. A comparison of our algorithms with convexity-constrained least squares based methods (e.g., \citet{deng2020inference}) and is also given in the Supplementary Material.

\section{Discussion}
\label{sec:discussion}

In the present paper, we studied optimal estimation and inference for the minimizer and minimum of a convex function under a non-asymptotic local minimax framework. We show in the Supplementary Material \citep[Section \ref{subsec:connection_classical}]{CaiChenZhuSupplement} that the results in this paper can be readily used to establish the optimal rates of convergence over the convex smoothness classes under the classical minimax framework.

A key advantage of the non-asymptotic local minimax  framework is its ability to characterize the difficulty of estimating individual functions
and to demonstrate novel phenomena that are not observable in the classical minimax theory. The Uncertainty Principle established in this paper highlights the fundamental tension between the accuracy of estimating the minimizer and that of estimating the minimum of a convex function. Similar results also apply to inference accuracy. It would be of great interest to establish uncertainty principles for other statistical problems, such as stochastic optimization with bandit feedback under shape constraints.

A more conventional approach to assessing the difficulty at individual functions is to consider the minimax risk over a local neighborhood of a given function in the parameter space. However, this entails defining a topology/metric over the parameter space and specifying the size of the local neighborhood, which can be challenging due to variations across problems.  Metrics such as $L_1$, $L_2$, or weighted $L_2$ distances are often employed. In contrast, our framework does not require specifying topology/metric on the parameter space or the size of a local neighborhood. This makes it more convenient and easier to use.  

The correct form of the conventional local minimax benchmark for the minimizer under \eqref{white.noise.model} is given by $R_{con,z}(\varepsilon;f) = \inf_{\hat Z}\sup_{ h\in  B(f, \varepsilon)}\E_{h} |\hat Z -Z(h)|$, where $B(f, \varepsilon) = \{h\in \FC: \|h-f\|_2\le \varepsilon\}$ is the $\varepsilon$-neighborhood of $f$ in $\FC$. The order of $R_{con,z}(\varepsilon;f)$ is clearly no smaller than our local minimax benchmark $R_z(\varepsilon;f)$ because
\beas
R_{con,z}(\varepsilon;f) &=& \inf_{\hat Z}\sup_{ h\in B(f, \varepsilon)} \E_{h} |\hat Z -Z(h)| \ge \inf_{\hat Z}\sup_{ h\in B(f, \varepsilon)}\sup_{g\in\{f,h\}}\E_{g} |\hat Z -Z(g)| \\
&\ge& \Phi(-0.5) \omega(\varepsilon;f) \sim R_z(\varepsilon;f),
\eeas
where the last inequality holds because of $\lim_{\eta\to \varepsilon} \sup_{g: \|g-f\|_2 = \eta } |Z(g) - Z(f)| = \omega(\varepsilon;f)$ and the last step is due to Theorem~\ref{thm:lowerbounds}. Analogous results hold for the other three problems. This means that our  benchmarks are at least as stringent as their conventional counterparts.  
Therefore, our procedures are also adaptive and optimal under the conventional local minimax framework.

The present work can be extended in different directions in addition to the aforementioned ones. For estimation,
the results can be easily generalized to the $\ell_q$ loss for $q > 1$.
It is interesting to consider the extremum under more general shape constraints such as $s$-convexity. In addition, estimation and inference for other nonlinear functionals such as the quadratic functional, entropies, and divergences under a non-asymptotic local minimax  framework can be studied.  We expect the penalty-of-superefficiency property to hold in these problems and our approach to be particularly helpful for the construction of the confidence intervals. Another important direction is to apply our  framework to other statistical models such as estimation and inference for the mode and the maximum of a log concave density function based on i.i.d. observations. We expect similar uncertainty principles to hold.

\section*{Acknowledgments}
We would like to thank the Associate Editor and the referees for their detailed and constructive comments which have helped to improve the presentation of the paper.

\begin{supplement}
\sname{Supplement to ``Estimation and Inference for Minimizer and Minimum of Convex Functions: Optimality, Adaptivity, and Uncertainty Principles"} 
\slink[url]{DOI:...}
	\sdescription{Section A of the supplement presents the simulation results. Section B discusses the comparisons with CLS based methods and the connection with the classical minimax framework. 
	The proofs of the main results and the technical lemmas are given respectively in Sections C and D.} 
\end{supplement}

\bibliographystyle{apalike}
\bibliography{Reference}

\begin{thebibliography}{}

\bibitem[Agarwal et~al., 2011]{agarwal2011stochastic}
Agarwal, A., Foster, D.~P., Hsu, D.~J., Kakade, S.~M., and Rakhlin, A. (2011).
\newblock Stochastic convex optimization with bandit feedback.
\newblock In Shawe-Taylor, J., Zemel, R., Bartlett, P., Pereira, F., and
  Weinberger, K.~Q., editors, {\em Advances in Neural Information Processing
  Systems}, volume~24. Curran Associates, Inc.

\bibitem[Auer et~al., 2007]{Auer2007ImprovedRF}
Auer, P., Ortner, R., and Szepesv{\'a}ri, C. (2007).
\newblock Improved rates for the stochastic continuum-armed bandit problem.
\newblock In Bshouty, N.~H. and Gentile, C., editors, {\em Learning Theory},
  pages 454--468, Berlin, Heidelberg. Springer Berlin Heidelberg.

\bibitem[Belitser et~al., 2012]{belitser2012optimal}
Belitser, E., Ghosal, S., and van Zanten, H. (2012).
\newblock Optimal two-stage procedures for estimating location and size of the
  maximum of a multivariate regression function.
\newblock {\em The Annals of Statistics}, 40(6):2850--2876.

\bibitem[Birge, 1989]{birge1989grenader}
Birge, L. (1989).
\newblock {The Grenader estimator: A nonasymptotic approach}.
\newblock {\em The Annals of Statistics}, 17(4):1532--1549.

\bibitem[Blum, 1954]{blum1954multidimensional}
Blum, J.~R. (1954).
\newblock {Multidimensional stochastic approximation methods}.
\newblock {\em The Annals of Mathematical Statistics}, 25(4):737--744.

\bibitem[Cai et~al., 2023]{CaiChenZhuSupplement}
Cai, T.~T., Chen, R., and Zhu, Y. (2023).
\newblock Supplement to ``{E}stimation and {I}nference for {M}inimizer and
  {M}inimum of {C}onvex {F}unctions: {O}ptimality, {A}daptivity, and
  {U}ncertainty {P}rinciples".

\bibitem[Cai and Low, 2015]{cai2015framework}
Cai, T.~T. and Low, M.~G. (2015).
\newblock A framework for estimation of convex functions.
\newblock {\em Statistica Sinica}, 25(2):423--456.

\bibitem[Cai et~al., 2013]{cai2013adaptive}
Cai, T.~T., Low, M.~G., and Xia, Y. (2013).
\newblock Adaptive confidence intervals for regression functions under shape
  constraints.
\newblock {\em The Annals of Statistics}, 41(2):722--750.

\bibitem[Chatterjee et~al., 2016]{zhu2016local}
Chatterjee, S., Duchi, J.~C., Lafferty, J., and Zhu, Y. (2016).
\newblock Local minimax complexity of stochastic convex optimization.
\newblock In Lee, D., Sugiyama, M., Luxburg, U., Guyon, I., and Garnett, R.,
  editors, {\em Advances in Neural Information Processing Systems}, volume~29.
  Curran Associates, Inc.

\bibitem[Chen, 1988]{chen1988lower}
Chen, H. (1988).
\newblock {Lower rate of convergence for locating a maximum of a function}.
\newblock {\em The Annals of Statistics}, 16(3):1330--1334.

\bibitem[Chen et~al., 1996]{chen1996estimation}
Chen, H., Huang, M.-N.~L., and Huang, W.-J. (1996).
\newblock Estimation of the location of the maximum of a regression function
  using extreme order statistics.
\newblock {\em Journal of Multivariate Analysis}, 57(2):191--214.

\bibitem[Deng et~al., 2020]{deng2020inference}
Deng, H., Han, Q., and Sen, B. (2020).
\newblock Inference for local parameters in convexity constrained models.
\newblock {\em arXiv preprint arXiv:2006.10264}.

\bibitem[Dippon, 2003]{dippon2003accelerated}
Dippon, J. (2003).
\newblock Accelerated randomized stochastic optimization.
\newblock {\em The Annals of Statistics}, 31(4):1260--1281.

\bibitem[Dumbgen, 1998]{dumbgen1998new}
Dumbgen, L. (1998).
\newblock New goodness-of-fit tests and their application to nonparametric
  confidence sets.
\newblock {\em The Annals of Statistics}, 26(1):288--314.

\bibitem[Facer and M{\"u}ller, 2003]{facer2003nonparametric}
Facer, M.~R. and M{\"u}ller, H.-G. (2003).
\newblock Nonparametric estimation of the location of a maximum in a response
  surface.
\newblock {\em Journal of Multivariate Analysis}, 87(1):191--217.

\bibitem[Ghosal and Sen, 2017]{ghosal2017univariate}
Ghosal, P. and Sen, B. (2017).
\newblock On univariate convex regression.
\newblock {\em Sankhya A}, 79(2):215--253.

\bibitem[Guntuboyina and Sen, 2018]{guntuboyina2018nonparametric}
Guntuboyina, A. and Sen, B. (2018).
\newblock Nonparametric shape-restricted regression.
\newblock {\em Statistical Science}, 33(4):568--594.

\bibitem[Hengartner and Stark, 1995]{hengartner1995finite}
Hengartner, N.~W. and Stark, P.~B. (1995).
\newblock {Finite-sample confidence envelopes for shape-restricted densities}.
\newblock {\em The Annals of Statistics}, 23(2):525--550.

\bibitem[Kiefer, 1982]{kiefer1982optimum}
Kiefer, J. (1982).
\newblock Optimum rates for non-parametric density and regression estimates
  under order restrictions.
\newblock {\em Statistics and Probability: Essays in honor of CR Rao}, 419:428.

\bibitem[Kiefer and Wolfowitz, 1952]{kiefer1952stochastic}
Kiefer, J. and Wolfowitz, J. (1952).
\newblock Stochastic estimation of the maximum of a regression function.
\newblock {\em The Annals of Mathematical Statistics}, 23(3):462--466.

\bibitem[Kleinberg, 2004]{kleinberg2004nearly}
Kleinberg, R. (2004).
\newblock Nearly tight bounds for the continuum-armed bandit problem.
\newblock In {\em Proceedings of the 17th International Conference on Neural
  Information Processing Systems}, NIPS'04, pages 697--–704, Cambridge, MA,
  USA. MIT Press.

\bibitem[Kleinberg et~al., 2019]{Kleinberg2019}
Kleinberg, R., Slivkins, A., and Upfal, E. (2019).
\newblock Bandits and experts in metric spaces.
\newblock {\em J. ACM}, 66(4).

\bibitem[Mokkadem and Pelletier, 2007]{mokkadem2007companion}
Mokkadem, A. and Pelletier, M. (2007).
\newblock A companion for the {K}iefer--{W}olfowitz--{B}lum stochastic
  approximation algorithm.
\newblock {\em The Annals of Statistics}, 35(4):1749--1772.

\bibitem[Muller, 1989]{muller1989adaptive}
Muller, H.-G. (1989).
\newblock {Adaptive nonparametric peak estimation}.
\newblock {\em The Annals of Statistics}, 17(3):1053 -- 1069.

\bibitem[Polyak and Tsybakov, 1990]{polyak1990optimal}
Polyak, B.~T. and Tsybakov, A.~B. (1990).
\newblock Optimal order of accuracy of search algorithms in stochastic
  optimization.
\newblock {\em Problemy Peredachi Informatsii}, 26(2):45--53.

\bibitem[Shoung and Zhang, 2001]{shoung2001least}
Shoung, J.-M. and Zhang, C.-H. (2001).
\newblock Least squares estimators of the mode of a unimodal regression
  function.
\newblock {\em The Annals of Statistics}, 29(3):648--665.

\end{thebibliography}

\end{document}